**Adaptive (re)operations facilitate environmental flow maintenance downstream of multi-purpose reservoirs**

Akshay Sunil[1], Riddhi Singh[1,2], Manvitha Molakala[1]

[1]Department of Civil Engineering, Indian Institute of Technology Bombay, Powai, Mumbai, Maharashtra 400076, India.

[2]Interdisciplinary Programme in Climate Studies, Indian Institute of Technology Bombay, Powai, Mumbai, Maharashtra 400076, India.

[*] Corresponding author:

E-mail address: akshaysunil172@gmail.com

# Adaptive (re)operations facilitate environmental flow maintenance downstream of multi-purpose reservoirs

**Abstract**

Multi-purpose reservoirs support socioeconomic development by providing irrigation, domestic water supply, hydropower, and other services. However, impoundment of water impacts instream aquatic ecosystems. Thus, the concept of minimum environmental flows (MEFs) was established to restore the benefits of naturally flowing rivers by specifying minimum flow rates to be maintained downstream of dams. But varying legislative contexts under which multi-purpose reservoirs operate may not always necessitate MEF releases. To what extent the release of MEF affects other sectoral benefits remains an open-ended and possibly a site-specific inquiry. A related issue is - how does the order in which releases are prioritized influences sectoral performances? We analyse these issues for the Nagarjuna Sagar reservoir, one of the largest multipurpose reservoirs in southern India. We formulate two versions of a multi-objective decision problem. PF_MEF formulation prioritizes MEF releases over releases for water demand satisfaction, followed by hydropower releases. PF_nMEF formulation follows the regional legislative rule releasing first for demand satisfaction, followed by hydropower and MEF releases. The objectives include: annual demand deficits (minimize), hydropower production (maximize), reliability of maintaining MEF (maximize), and reliability of non-exceedance of high flows downstream (maximize). A dynamic rule using radial basis functions (RBFs) specifies releases as a function of reservoir storage and the Borg multi-objective evolutionary algorithm is used to identify the Pareto approximate solutions. The reliability of MEF satisfaction, annual demand deficits, and hydropower production attained across the Pareto-optimal strategies obtained for PF_MEF (PF_nMEF) formulation is 86–99 % (63 - 80%), 47-1063 (17.6- 340) Mm$^3$ , and 3452-3703 (3408-3650) GWh, respectively. Results thus indicate that prioritizing MEF releases improves can meet MEF requirements without significant compromises in other objectives. We hypothesize that similar investigations may reveal how simple modification of release order may improve ability of other reservoirs to meet environmental goals.

Keywords: water resources management, minimum environmental flow, MEF, multi-purpose reservoir, rule curve, multi-objective optimization

1. Introduction

Global water requirements for food production, power generation, domestic water supply, and industrial needs have been growing over the last few decades (Hanjra and Qureshi, 2010). This has increased investments in large scale water infrastructure such as multi-purpose reservoirs that support domestic water supply, irrigation, hydropower production, navigation, recreation, fisheries, etc. However, the construction of dams or barrages presents an obstruction to the natural flow of rivers, putting ecological water needs at risk due to their hindrance of the river continuum (Dunne and Leopold, 1978; Vannote et al., 1980; NRC, 1992; Poff et al., 1997; Nilsson et al., 2005; Loures et al., 2012; Hairan et al., 2021; Mezger et al., 2021; Hoque et al., 2022). "Large hydropower dams alter the natural flow of rivers by producing a rapidly fluctuating hydropower regime downstream when operated for "peaking power" (Watts et al., 2011). Similarly, water supply dams alter the seasonality of river flows by extending low flow periods during the dam filling phase and introducing unnaturally high flows during the release phase (Watts et al., 2011)".

The concept of minimum environmental flows (MEF) arose in response to this conflict between human and ecological needs, providing guidance to reservoir operators on the minimal flow rates required downstream of major structures to sustain riverine ecology. The widely cited Brisbane Declaration (2007) defines environmental flows as the quantity, quality, and timing of water, sediments, and biota flows required to sustain freshwater and estuarine ecosystems, as well as human livelihood and wellbeing. MEF releases from reservoirs is essential in downstream channel maintenance, upkeeping the riparian vegetation, bird breeding, and wetland flooding (Naiman et al., 1995; Poff et al., 1997; Naiman et al., 2007; Dyer et al., 2014). Preserving the temporal variability of river flows is crucial in supporting a range of ecosystem processes that sustains native biodiversity (Poff and Ward, 1989; Richter et al., 1996; Poff et al., 1997; Postel et al., 1997; Puckridge et al., 1998; Richter, 2010; Hayes et al., 2018). Maintaining MEFs also enables ecosystem services such as provisioning, regulating, cultural, and supporting services, benefiting human society (Clay, 1995; Forslund et al., 2009; Larinier et al., 2002). Numerous studies show that restoring the components of the river's natural flow by dam reoperation can reverse the adverse effects caused by dams (Nyatsanza et al., 2015; Opperman et al., 2019; Dutta et al., 2020; Merritt et al. 2009; Poff and Zimmerman, 2010; Ahn et al., 2018; Owusu et al., 2020, Owusu et al., 2021). Owusu et al. (2020) identified 50 river systems with 69 dams which were re-operated for releasing MEFs and concluded that MEF requirements are not always conflicting with other objectives such as water supply, hydropower production, and flood protection. A few studies have quantitatively explored the trade-offs between hydropower, water supply, and environmental flow protection in dam operation (Sandoval-Solis and McKinney, 2014; Uen et al., 2018; Al-Jawad et al., 2019; Jiang et al.2019; Quinn et al.2017; Quinn et al., 2019; Derepasko et al., 2021; Veena et al. 2021). Some studies model MEF requirements as minimum flow constraints rather than as objective functions (Sale et al., 1982; Castelletti et al., 2008; Yin et al., 2010; Desreumaux et al., 2018; Fatahi Nafchi et al., 2021; Naghdi et al., 2021). However, releases for satisfying water supply or irrigation demand have generally been given a higher priority over other objectives (Sandoval-Solis and McKinney, 2014; Feng et al., 2018; Uen et al., 2018; Al-Jawad et al., 2019; Zhang et al., 2019; Naghdi et al., 2021, Jordan et al., 2022). For example, Al-Jawad et al. (2019) analysed the multi-sectoral trade-offs between demand satisfaction, hydropower production, and MEF releases for the Himren dam reservoir in the Diyala River basin, Iraq. In their analysis, these multi-sectoral objectives were simultaneously optimized using two multi-objective evolutionary algorithms, the Borg MOEA and ε-DSEA. Demand satisfaction was prioritized, followed by hydropower generation and environmental flow releases. Similarly, Uen et al. (2018) used the Non-dominated sorting genetic algorithm II (NSGA-II) to optimize daily reservoir operations of the Shihmen reservoir on the Tamsui River, Taiwan. Here too, the releases for domestic demand satisfaction were according to the highest priority followed by releases for irrigation demand satisfaction and hydropower production. Therefore, MEF releases have been generally accorded lower priority among multi-sectoral releases in analytical studies on reservoir performance. However, it has not yet been explored whether the priority of release decisions alters the nature of trade-offs in a water resources optimization problem involving multi-sectoral water uses.

Recent studies call for reservoir reoperations across the globe to restore ecosystems (Watts et al., 2011; Yang and Cai, 2011; Nyatsanza et al., 2015; Rheinheimer et al., 2016; Chen and Olden, 2017; Opperman et al., 2019). As already shown by Owusu et al. (2020), the common assumption that there is a conflict between social and environmental objectives may not hold across different reservoir settings. Furthermore, the magnitude of the conflict may not be as perceived and it may also be feasible to reduce the conflict by realtering dam operation strategies. Warner et al. (2014) reviewed the dams where environmental flow maintenance was implemented through adaptive reservoir re-operations as a part of the Sustainable Rivers Project (SRP), which began in 2002. This study concluded that there exist opportunities to modify reservoir operations to improve environmental conditions as well as traditional reservoir purposes such as water supply and hydropower generation. Recently, a few studies have quantified the multi-objective trade-offs for run-of-the-river plants where there are minimal flow alterations (Basso and Botter, 2012; Kuriqi et al., 2019; Owusu et al., 2023). Basso et al. (2023) quantified the trade-offs between hydropower generation and environmental flow satisfaction for a mini-hydro plant in the Piave River catchment of northern Italy. Their findings highlight the critical role of streamflow variability in deciding the resultant multi-objective trade-offs. Kuriqi et al. (2019) quantified the trade-offs between ecological and hydropower objectives for different environmental flow scenarios for 20 run-of-the-river plants within five river basins in the Iberian Peninsula. They generated eight scenarios of environmental releases using different hydrology-based environmental flow methods. They concluded that a balance can be attained between hydropower generation and environmental flow satisfaction for certain environmental flow scenarios. Thus, recent literature points towards a likely benefit in quantifying the trade-offs involved in re-operating multi-purpose storage reservoirs for prioritizing MEFs, even though there has never been a direct comparison of resultant trade-offs between cases that prioritize and do not prioritize MEF."

Despite their increasingly important role in maintaining riverine ecology, there are no globally accepted recommendations on how to (re)operate multi-purpose reservoirs to meet MEF goals alongside other objectives (Chen et al., 2013). The Brisbane Declaration (2007) emphasized that governments should put environmental flow standards in place and enforce them. Several other international NGOs including the International Union for Conservation of Nature (IUCN) also embraced the concept and stressed its relevance in emerging economies (Dyson et al., 2008). In response to this, Australia, the European Union, South Africa, and the USA have integrated environmental flow provisions into their legislation and policy-making frameworks (Venot et al., 2008). The development of methodologies for MEF assessment in the United States began in the late 1940s to incorporate new environmental and freshwater legislation, and now one-third of the state's rivers are designated for special protection

(Krchnak et al. 2009; Poff and Matthews, 2013). In Australia, the Water Act (2007) established diversion limits on river flows for agricultural and other consumptive purposes (Water Resources Environmental Flow Guidelines, 2019). The European Union has approved the Water Framework Directive (WFD) to this end (Souchon and Keith, 2001, Acreman and Ferguson, 2010, Ramos et al., 2018). Similarly, in South Africa, water is set aside as a reserve with its own legal standing under the National Water Act (NWA). It is evident that there is a significant variation in the prescribed approaches for protecting MEFs, while little guidelines exist for water managers on how to operate the reservoir considering other needs.

In this study, we investigate how prioritization of the order of release decisions affects multi-objective trade-offs among competing objectives in reservoir operations, specifically focusing on their ability to maintain MEF requirements. We use the Nagarjuna Sagar (NS) reservoir on the Krishna River basin in Southern India as an example of a typical multi-purpose reservoir that supplies water for irrigation as well as domestic needs while also generating hydropower. The NS reservoir is representative of large-scale projects with multiple objectives facing significant hydroclimatic variability. In recent years, the reservoir has faced several multi-year droughts during which these objectives have come into conflict leading to socio-political issues (Kumar and Jayakumar, 2020; Thomas et al., 2019). In addition, the reservoir is likely to participate in major inter-basin water transfers, that will further affect MEF releases (Veena et al. 2021; Molakala and Singh, 2022). Therefore, we investigate whether changing the prioritization order of releases may play a role in improving the ability of the NS reservoir in maintaining MEF downstream. For the NS reservoir, we answer the following questions:

1. What, if any, are the trade-offs between the main objectives of the NS reservoir?
2. Does the prioritization of MEF in release operations significantly affect the performance of other objectives?
3. How do varying definitions of MEF affect the answers to 1-2 above?

## 2. Study area and datasets

### 2.1 The Nagarjuna Sagar multi-purpose project

Due to the highly seasonal nature of the Indian summer monsoon, numerous reservoirs have been constructed in India to effectively use water. The Nagarjuna Sagar (NS) dam was built across the Krishna River in 1967 to bolster irrigation and hydroelectricity production in the Krishna River basin (Fig. 1). The Krishna River is India's fourth-largest river, and the Krishna Basin is the country's fifth-largest basin with a catchment area of 2,20,705 $Km^2$ at the dam site. The project consists of a dam with a live storage capacity of 5733 $Mm^3$ and two main canals, the right main canal (NSRC) and left main canal (NSLC), each with a hydropower capacity of 60 and 90 MW, respectively. The main power plant downstream of the dam has a hydropower potential of 960 MW (Karnatapu et al., 2020). The releases into the canals for irrigation water supply are also used for generating hydropower. The project has a total command area of 8955 $Km^2$ and nearly 87% (7435 $Mm^3$) of the water from the reservoir is used for irrigation while the remaining is used to satisfy domestic and industrial demands (Veena et al., 2021; Veena et al., 2022). Irrigation demands peak during the monsoon season when water-intensive crops are grown, then fall and reach a minimum during the summer. Since 2004, the reservoir also supplies nearly 40% of the total domestic water demand of the pharmaceutical hub of Hyderabad, located 114 km away (George et al., 2009).

Prior to 1982, the average annual inflow at the reservoir site was about 29,048 **$Mm^3$**, when another major dam was constructed upstream of the NS dam (Kumar and Jayakumar, 2020). Thus**,** around 30% of the inflows were to be used for irrigation and other purposes as per the original plans, leaving enough quantity for environmental flow requirements downstream. However, in recent years, a significant reduction in inflows to the reservoir due to upstream diversions has resulted in a further increase in the fraction of water diverted for human use at the dam site (Biggs et al., 2007; Venot et al., 2008; Venot, 2009; Gaur et al., 2007; Veena et al., 2021). While the NS dam has facilitated socioeconomic growth in the region surrounding it, there has also been a significant impact on the biodiversity of the Krishna River basin due to these developments (Smakhtin et al., 2007). Climate change is likely to further impact inflows to the NS reservoir, resulting in uncertain impacts on environmental flow-related releases (Maddu et al., 2022). Therefore, understanding the trade-offs between the socioeconomic and ecological objectives of the NS reservoir would be crucial to aid NS dam operations in the future.

**"Insert Figure 1"**

Fig. 1. (a) Location map of the Nagarjuna Sagar dam and its command area in India. (b) The land use in the Nagarjuna Sagar reservoir command area as of year 2016. Land use data was obtained from NRSC (National Remote Sensing Centre, https://bhuvan-app1.nrsc.gov.in/thematic/thematic/index.php)

**2.2 Data sources**

We obtained the inflow data for the NS reservoir for the period 1967 to 2016 from the Telangana Irrigation and Command Area Development Department (Veena et al., 2021; Molakala and Singh, 2022; Veena et al., 2022). Monthly irrigation demands were quantified by calculating the net irrigation requirements for each crop assuming a fixed annual cropping pattern (Prasad et al., 2006; Veena et al., 2021; Molakala and Singh, 2022). The annual domestic demands for the NS reservoir were calculated using population data and per capita demand estimates for Hyderabad (NWDA (National Water Development Agency) 2020) and were estimated to be 1000 **$Mm^3$** (Veena et al., 2021). The thresholds for the minimum environmental flow requirements were set at 30% of the mean annual flows (MAF) (Smakhtin and Anputhas, 2006; Karimi et al., 2021; Veena et al., 2021; Molakala and Singh, 2022). In addition, high flow exceedance thresholds downstream of the reservoir are defined by calculating the maximum daily downstream releases as estimated from historical downstream release data of the reservoir (Veena et al., 2021; Molakala and Singh, 2022).

## 3. Methods

In complex decision contexts such as planning large-scale water infrastructure, exploring alternative framings of a problem enable decision makers to analyse the consequences of their actions across differing assumptions of models, objectives, and/or inputs (Quinn et al., 2017). Through the development of alternative quantitative representations of objective functions, the effect of competing problem formulations on trade-offs in the management of water resources of the Red River basin in Vietnam was assessed by Quinn et al. (2017). Since then, several studies have adopted this approach to analyse water resource decision problems. For example, Veena et al. (2021) adopted four framings of an inter-basin water transfer problem to quantify the value of coordinated information exchange in improving the performance of water transfer strategies. Thus, we also framed alternative formulations of the reservoir release decision problem to understand whether the order of prioritization of releases from a multi-purpose reservoir affects the multi-sectoral performances across demand and environmental flow related objectives. Two framings of the reservoir release decision problem evaluate the effect of release prioritization decisions on trade-offs between demand satisfaction, hydropower generation, minimum environmental flow satisfaction, and flood control. The PF_MEF formulation prioritizes MEF releases over demand satisfaction while the PF_nMEF formulation does not. In both cases, order of other releases is set as that legislatively prescribed: demand releases followed by hydropower releases (Fig. 2).

Within each framing, we employ a dynamic reservoir release rule for hydropower production and identify (Pareto-)optimal operating strategies using a meta-heuristic algorithm. Section 3.1 describes the systems model and the multiple objectives of the decision problem. Section 3.2 provides a detailed discussion of the control policy formulation for the NS dam and the multi-objective optimization experiment.

**"Insert Figure 2"**

Fig. 2. Two alternative formulations of the systems model for the NS reservoir operations. (a) PF_nMEF prioritises demand releases over minimum environmental flows, and (b) PF_MEF prioritises environmental flows over demands and other releases.

### 3.1 Systems model and reservoir management objectives

A dynamic water balance model tracks the reservoir storage at each timestep after accounting for various releases and inflows (Equation 1).

$$s_{t,j} = s_{t-1,j} + q_{t,j} - efr_{t,j} - dr_{t,j} - hpr_{t,j} - ewr_{t,j} \quad (1)$$

In Equation (1), $s$ is the reservoir storage, $q$ is the inflow, $efr$ is the release for MEF, $dr$ is the release for satisfying demands, $hpr_t$ is the release for hydropower production, and $ewr$ is the excess release. Excess water is released downstream when water level in the reservoir exceeds the live storage capacity of reservoir. The subscript 't' denotes the time step of the model (daily), and 'j' denotes the synthetic inflow realization. The model is run using an ensemble of 10 years of daily streamflow time series synthetically generated from 35-year historically observed inflows from 1968-2003. A total of 10,000 synthetic inflow time series are generated following the Kirsch-Nowak streamflow generator (Nowak et al., 2010; Kirsch et al., 2013) applied to the NS reservoir inflow data by Veena et al. (2021) (See

Supplementary Text S1 for details). Using these ensembles allows us to explore a wide range of high and low flows, thus enabling a generalization of the results beyond the observed inflow records.

Based on the historical context of the NS dam, we formulated four main objectives of the NS reservoir operations: annual average hydropower generation (maximize), annual average demand deficit (minimized), reliability of maintaining MEFs (maximize), and reliability of avoiding floods downstream of the NS dam (maximize). The NS dam was constructed with the primary goal of irrigation and domestic water supply followed by hydropower production. However, given the context of our analysis, we included an environmental flow related objective to understand the consequences of dam operations on MEF releases. Furthermore, the Krishna River experienced three major floods between the period 1969-2013 (Killada et al., 2012; Veena et al., 2021). As flood protection requires maintaining empty reservoir conditions to absorb flood waters, while hydropower production as well as demand satisfaction requires maintaining a high head (and thus storage), these may come in conflict. We therefore included a flood related objective in the analysis. The mathematical formulation of each objective is listed in Table 1. These objectives are aggregated across the inflow ensembles by estimating the mean value of the objective function.

Table 1. The objective functions considered for the NS project.

| SNo | Objective | Equation | Notation |
|---|---|---|---|
| 1 | Annual average hydropower production ($J_{hydropower}$) [GWh] | $J_{hydropower} = \frac{1}{NR}\sum_{j=1}^{NR}\frac{1}{m}\sum_{t=1}^{T} N_{t,j}$  $N_{t,j} = min(\eta\, \gamma\, hpr_{t,j}\, H_{t,j},\, HP_{cap})$ | $NR$: number of synthetic inflow realizations  $m$: number of years in the in the planning horizon |
| | | | $N_{t,j}$: power output (GW) at time $t$ for realization $j$ |
| | | | $T$: Total number of time steps |
| | | | $\eta$: turbine efficiency |
| | | | $\gamma$: specific weight of water |
| | | | $H_{t,j}$: net reservoir head at time $t$ for realization $j$ |
| | | | $HP_{cap}$: maximum hydropower capacity |
| 2 | Annual average demand deficit ($J_{avg\ deficit}$) [Mm$^3$] | $J_{avg\ deficit}$ $= \frac{1}{NR}\sum_{j=1}^{NR}\frac{1}{m}\sum_{t=1}^{T} dr_{t,j} - DD_{t,j}$ | $DD_{t,j}$: total water demand at time $t$ for realization $j$ |
| 3 | Reliability of minimum environmental flow ($J_{EF}$) [%] | $J_{EF} = \frac{1}{NR}$ $\sum_{j=1}^{NR}\frac{1}{T}\sum_{t=1}^{T} EF_{t,j} \times 100$ $EF_{t,j}$ $= \begin{cases} 1 \text{ if } (efr_{t,j} < mef_t) \\ 0 \quad\quad\quad else \end{cases}$ | $mef_t$ : minimum environmental flow threshold at time $t$ |
| 4 | Reliability of non-exceedance of high flows ($J_{FT}$) [%] | $J_{FT} =$ $\frac{1}{NR}\sum_{j=1}^{NR}\frac{1}{T}\sum_{t=1}^{T} FT_{t,j} \times 100$ $FT_t = \begin{cases} 1 \text{ if } (r_{t,j} < ft_t) \\ 0 \quad\quad\quad else \end{cases}$ | $r_{t,j} = efr_{t,j} + hpr_{t,j} + ewr_{t,j}$ : actual downstream release at time $t$ for realization $j$ $ft_t$ : high flow threshold at time $t$ |

## 3.2 Control policy formulation and multi-objective optimization

Previous research has shown that radial basis function (RBF) based dynamic adaptive rule curves that condition release decisions on storage states perform well for reservoir operation problems (Salazar et al., 2016; Quinn et al., 2019). The strong universal approximation mathematical properties of Gaussian RBFs have shown to be a very effective in generating flexible operational rules (Giuliani et al., 2014; Giuliani et al.,2016). Giuliani et al. (2014) showed that RBF parameterization of releases outperforms Artificial Neural Networks (ANN) in exploring the entire trade-off space between conflicting objectives. We thus use RBFs to condition the decision of hydropower releases on the storage state of the NS reservoir (Equation 3).

$$\boldsymbol{hpr_{t,j}} = \boldsymbol{Q_{HP}}\, exp\left[-\frac{(ns_{t,j}-c)^2}{r^2}\right] \qquad (2)$$

In Equation 2, $\boldsymbol{Q_{HP}}$ is the maximum capacity for hydropower releases, $\boldsymbol{ns_{t,j}}$ is the normalized storage, $\boldsymbol{c}$ and r are the centres and radii of the RBF, and other terms have been described previously. c and r are the two parameters to be optimized and their ranges are [-1,1], and [0,1] respectively. The normalized storage is obtained by dividing the storage for a timestep $t$ by the sum of live storage and the maximum inflow across the realizations.

We then set up a multi-objective optimization experiment to identify (Pareto-)optimal sets of RBF parameters (Equation 3).

$$\boldsymbol{\theta}^* = \arg\min_{\boldsymbol{\theta}} J(\boldsymbol{\theta}) \qquad (3)$$

In Equation 3, $\boldsymbol{\theta}$ is the decision variable vector comprising of various possible combinations of b and c, J is the vector of objective function values obtained for a given $\boldsymbol{\theta}$. The non-dominated set of decision variable vectors $\boldsymbol{\theta}^*$ that minimizes $J(\theta)$ constitute the Pareto-optimal release strategies. The Borg multi-objective evolutionary algorithm (MOEA) is used to search the decision variable space to identify Pareto-optimal water release strategies for this multi-objective problem (Hadka and Reed, 2013). The Borg-MOEA has six auto-adaptive recombination operators and uses the ε-box dominance concept for dominance sorting (Kollat and Reed, 2006). It has been shown to outperform other MOEAs on a variety of water resource problems (Hadka and Reed, 2012; Hadka et al., 2012; Hadka and Reed, 2013; Salazar et al., 2016).

We ran the optimization experiment for 10,000 function evaluations (NFE) using epsilon values for the hydropower production (GWh), average annual deficit (Mm$^3$), MEF reliability [%], and reliability of non-exceedance of high flows [%] as (50,50,0.01,0.01), respectively. To ensure convergence of the solutions, we ran the analysis across five random seeds and monitored the hypervolume metric (Supplementary Material S2). In a single evaluation of the objective function vector, we employ a stochastic search procedure that involves randomly selecting 100 inflow time series out of a total of 10,000 and evaluating the mean value of the objective function across these. This approach considerably reduces the computational costs compared to aggregating inflows across all 10,000 realizations in each function call of the algorithm (Singh et al., 2015; Veena et al., 2021).

### 3.3 Sensitivity of results to the choice of MEF thresholds

In order to test for the sensitivity of the results to the choice of MEF thresholds, we varied the percentage cut-off requirement prescribed by Smakhtin and Anputhas (2006) from 30% to 80% of MAF, and re-optimized both formulations for the new requirements. Additionally, we determined the MEF thresholds based on the standards set by a central legislative body (MoEF&CC, 2006, 2020). For computing the MoEF&CC thresholds, a 35-year daily historical inflow series from 1968 to 2003 for the NS reservoir was used. First, the total annual flow is estimated for each year and then using the Weibull's formula, the total annual flow value with 90% exceedance probability is identified. Then, prescribed thresholds for monsoon (30% for June to September) lean (25% for October, April, and May) and dry (20% for November to March) months are applied to the daily data for the year corresponding to this total annual flow (CWC, 2019). All computations are carried out for water year from June to May.

## 4. Results

To understand the effect of priority order of releases on MEF satisfaction, we first compare the objective function performance and resultant reservoir operation strategies across the two formulations, one that prioritizes MEF (PF_MEF) and another that does not (PF_nMEF, Section 4.1). Here, we compare the multi-sectoral trade-offs between hydropower generation, average demand deficit, flood reliability, and MEF reliability across the two alternative formulations. We quantify and visualise the objective trade-offs to analyze how much additional demand satisfaction or hydropower production loss results from prioritizing MEF. We then explore the storage-release dynamics of reservoir reoperation to explain the variation of performance across the strategies obtained from PF_MEF and PF_nMEF formulations (Section 4.2). To this end, we select example compromise strategies from each problem formulation (PFMEF_S1, and PFnMEF_S2) to show the dynamics of operating the reservoirs across three climate regimes (dry, normal, and wet years). Furthermore, to test for the sensitivity of the results to the choice of MEF thresholds, we re-optimize the PF_MEF formulation for varying MEF requirements (Section 4.3). We then visualise how the multi-objective trade-offs are affected across these varying thresholds and identify at which level of MEF requirement it becomes challenging to achieve a balance between these trade-offs.

### 4.1 Multi-sectoral trade-offs for PF_MEF and PF_nMEF

On performing stochastic multi-objective optimization for the PF_MEF and PF_nMEF formulations, we obtained a total of 30 Pareto-optimal water release strategies, 14 for the PF_MEF formulation (green) and 16 for the PF_nMEF (red) formulation (Fig. 3). These strategies were re-evaluated against a larger ensemble of 100,000 inflow realizations to test the performance against a wider sample of historical inflow variability than that used during the optimization. The mean value of the objective functions across these larger ensembles are presented here. We also Pareto-sorted these solutions, taking into consideration all objectives and identified 6 non-dominated solutions from each formulation (solid solutions). To visualize performance of these strategies across multiple objectives, we use the parallel axes plot (Fig. 3a), pair-wise performance plots (Figs. 3b-e) and the cumulative distribution function (CDF) of each objective function across each formulation (Figs. 3f-i). Each vertical axis of the parallel axis represents an objective, with the preferred direction being upward. Each line connecting across the vertical axes represents a Pareto-optimal strategy. The points at which these lines intersect each vertical axis represent the performance of each strategy for that objective. A straight line intersecting all vertical axes at the top is an ideal solution.

We find that the Pareto optimal solutions from the PF_MEF formulation maintain reasonable performance across all four objectives while improving median MEF reliability to 93% across all re-sorted when compared to the PF_nMEF formulation (77% median MEF reliability, Figs. 3f-i). Solutions that prioritize MEF also perform better in the hydropower objective when compared to the solutions from the PF_nMEF (Fig. 3f). The median objective function values for hydropower production are 3662 (3457) GWh for the PF_MEF (PF_nMEF) across all solutions, respectively. The range of annual hydropower generation and MEF reliability objective values for the solutions from the PF_MEF (PF_nMEF) formulation is 3452- 3703 GWh (3408-3650 GWh), and 86–99% (63 –80 %). However, this improvement comes at a trade-off with the PF_MEF solutions performing worse in the average demand deficit objective. The median value of the average demand deficit attained across all solutions is 709 Mm$^3$ for the PF_MEF formulation, and 60 Mm$^3$ for the PF_nMEF formulation. Of note, however, is that these are summary statistics across all solutions. But there are unique solutions in the PF_MEF formulation which have reasonable performance in the demand deficit objective while maintaining better MEF performance, thus providing an alternative way to operate the NS reservoir.

Our findings show that adaptive operational strategies can be developed by prioritizing environmental flows without significantly risking other objectives. The strategies from the PF_nMEF formulation result in annual demand deficit ranging from 18-340 **Mm³**, while those from the PF_MEF formulation have a range of 47-1063 **Mm³**. Thus, the lowest demand deficit values are quite comparable for both formulations. Simultaneously, the strategies that attain the lowest demand deficits result in MEF reliability of 99% and 80% for the PF_MEF and PF_nMEF formulation, respectively. Note also that the reliability of non-exceedance of high flows objective attained a rather high value (>97%) across all solutions indicating that there are no significant issues related to flood management considering historical variation of inflows.

We also find a significant trade-off between hydropower and demand deficit objectives (Fig. 3b). For the PF_MEF formulation, the strategy that corresponds to the highest hydropower generation (3703 GWh) results in an average demand deficit of 1021 **Mm³** while the strategy with the least demand deficit (46.8 **Mm³**) results in lower hydropower generation of 3452 GWh. For the PF_nMEF formulation, the strategy that corresponds to the highest hydropower generation (3650 GWh) results in an average demand deficit of 340 **Mm³** while the strategy with the least demand deficit (17.6 **Mm³**) results in lower hydropower generation of 3408 GWh. Thus, 251 GWh hydroelectricity corresponds to 974.4 **Mm³** of water deficits in the PF_MEF formulation, while 243 GWh hydroelectricity corresponds to 323 **Mm³** of water deficits in the PF_nMEF formulation. Thus, while trade-offs exist between these objectives, their quantification does depend upon the specific problem formulation. This trade-off arises because when more water is supplied to meet demands, the reservoir's available head to generate electricity decreases resulting in lesser hydropower production. A similar trade-off is also visible between MEF reliability and hydropower generation. Note also that the objective function values calculated using the historical releases data with observed inflow for the NS dam operations from 1968-1983 are 1784 GWh, 1304 Mm³, 98%, and 67% for hydropower production, demand deficit, reliability of non-exceedance of high flows , and MEF reliability, respectively (Fig. 3, black markers). This indicates that dynamic adaptive operational strategies have the potential to improve the performance of the NS dam.

"Insert Figure 3"

Fig. 3. (a) Parallel axes plot illustrating the trade-offs between hydropower generation [GWh], average demand deficits [$Mm^3$], reliability of non-exceedance of high flows (flood reliability) [percentage], and MEF reliability [percentage] for the NS reservoir. Arrows represent the preferred direction for each objective function. Each line represents a solution obtained by multi-objective optimization. Solid lines represent the non-dominated solutions across the solutions obtained from the PF_MEF (green) and PF_nMEF (red) formulations. Two strategies from each formulation are further studied, identified by star markers. Dashed lines represent the objective function performance obtained by evaluating the starred strategies for the observed inflow data. Black line represents the objective function values determined for the historical period computed using the historical release data. (b-e) Pairwise plots between the four objectives. (f-i) The cumulative distribution function (y-axis) vs. objective function values (x-axis) across the non-dominated strategies sorted for all objectives as in panel a.

We now select one solution each from the two formulations for further analysis. These selected strategies are named PFMEF_S1, and PFnMEF_S2 and are from the PF_MEF and

PF_nMEF formulations respectively (stars in Fig. 3). These strategies are selected as they attain a comparable objective performance in hydropower generation, average deficit, and reliability of non-exceedance of high flows objectives with the PFMEF_S1 strategy being much improved in the MEF reliability objective. Thus, PFMEF_S1 is likely to be of interest to stakeholders by offering a very high improvement in maintaining MEFs and minimally degrading the annual deficit targets. These two selected solutions also attain the lowest annual average deficit for their respective formulations.

### 4.2 Reservoir storage-release dynamics for PF_MEF vs. PF_nMEF

We now examine the time variability of reservoir fluxes and subsequently the storage-release dynamics to explain the difference between the objective function performance of PFMEF_S1 and PFnMEF_S2 strategies, after re-evaluating them for historically observed inflows for 1968-1983. We first visualize the time series of hydropower releases, MEF releases, demand deficits, and reservoir storage of the selected strategies for representative dry (1972), normal (1979), and wet (1975) years during the observed period (Fig. 4). The year 1972 recorded the lowest value of annual inflow (15,859 $Mm^3$) while the year 1975 recorded the maximum value (74,805 $Mm^3$). The year 1979 is chosen as a typical normal year as it is the year with 50% exceedance probability of annual inflow (41,521 $Mm^3$). Reservoir related fluxes generally follow the intra seasonal variability of the Indian summer monsoon and associated water demand patterns. Even though both the selected strategies exhibit comparable performance for hydropower generation (3284 GWh for PFMEF_S1 and 3419 GWh for PFnMEF_S2), hydropower releases are consistently lower for the PFMEF_S1 strategy. The annual average hydropower release for the PFMEF_S1 strategy is 7629 $Mm^3$ compared to 14482 $Mm^3$ for the PFnMEF_S2 strategy. However, the PFMEF_S1 strategy maintains higher levels of water storage in the reservoir, with an average storage volume of 3361 $Mm^3$, compared to the PFnMEF_S2 strategy with an average storage of 2543 $Mm^3$ thereby preventing the storage level from dropping below a certain threshold for a majority of time periods. The high storage of PFMEF_S1 strategy ensures a higher water head which results in comparable hydropower generation even with lower release volumes. Note also that for these simulations based on observed inflows, the demand deficits are similar across both the strategies (162 $Mm^3$ for PFMEF_S1 strategy and 163 $Mm^3$ for PFnMEF_S2 strategy).

The high MEF reliability for the PFMEF_S1 strategy is attributed to the fact that MEF releases meet the MEF thresholds, especially during the monsoon months (June to September) when MEF requirements are higher (Fig. 4b). Note also that the PFMEF_S1 strategy was able to meet the MEF thresholds during the high MEF requirement monsoon months for all three years thus maintaining the intra-annual variations essential in river flow. On the other hand, PFnMEF_S2 strategy meets MEF requirements only during November-January for the dry year, September-May for the normal year, and August-May for the wet year. So, even in the wettest year on record, the PFnMEF_S2 strategy does not meet MEF requirements in the months of June and July. Owing to MEF prioritization, the reservoir's storage level for the PFMEF_S1 strategy drops below the PFnMEF_S2 strategy storage level in the dry year. This results in a higher demand deficit for the PFMEF_S1 strategy (2253 $Mm^3$) compared to the PFnMEF_S2 strategy (1123 $Mm^3$) for the dry year. In the normal and wet years, neither strategy results in demand deficits. Note also, the MEF thresholds used in this study for optimization are the thresholds suggested by Smakhtin and Anputhas (2006) which are considerably higher with an annual MEF requirement of 10164 **$Mm^3$** compared to the thresholds calculated using MoEF&CC standards with an annual MEF requirement of 2174 **$Mm^3$**. When considering the MOEF&CC standards (blue line in Fig. 4b), both strategies perform well in general except for the months of July-August when only the

PFMEF_S1 strategy meets the standards in all months while the PFnMEF_S2 strategy fails to meet the standards.

<div align="center">**"Insert Figure 4"**</div>

Fig. 4. (a) Inflows, (b) MEF releases, (c) hydropower releases, (d) demand deficits, and (e) reservoir storage (y-axis) of the NS dam against time (in days, x-axis) for PFnMEF_S2 (red colour) and PFMEF_S1 (green colour) strategies for a representative dry (panel 1), normal (panel 2), and wet year (panel 3). All variables are in volumetric units [$Mm^3$]. Star solutions from Fig. 3 (PFMEF_S1 and PFnMEF_S2) are re-evaluated for observed inflows from 1968-1983 for water years starting from June from which the representative dry (1972), normal (1979), and wet (1975) years are plotted. MEF thresholds computed using Smakhtin and Anputhas (2006) are shown by black line while MoEF&CC guidelines are shown by dark blue.

Upon investigating the temporal variability in MEF threshold satisfaction across the 100,000 inflow realizations, we find that MEF thresholds are satisfied for the PFMEF_S1 strategy in 97% of the realizations during the monsoon months (June to September). On the other hand, the PFnMEF_S2 strategy meets MEF requirements only in 41.6% of the realizations. We show this by plotting the reliability of meeting MEF threshold for each month across the stochastic inflow realizations (Fig. S3). During the non-monsoon months, the PFMEF_S1 strategy attains an average MEF reliability of 100%, while the PFnMEF_S2 strategy has an average MEF reliability of 97.3 %. This is due to overall low MEF requirements during the non-monsoon months.

The storage release dynamics for the PFMEF_S1 and PFnMEF_S2 strategies strongly suggest that a rather simple reoperation of the NS dam is likely to yield benefits for both hydropower and MEF related objectives (Fig. 5). The hydropower releases for both PFMEF_S1 and PFnMEF_S2 are determined by flexible rules owing to the nonlinear nature of radial basis functions. The hydropower releases monotonically increase as a function of the storage states, capped at the maximum daily discharge capacity of 101.9 $Mm^3$. For the PFMEF_S1 strategy, water is released for hydropower only when a storage head of 2000 $Mm^3$ is reached. In this way, the PFMEF_S1 strategy always maintains a sufficient storage head is the reservoir ensuring hydropower production along with release of MEF. On the other hand, PFnMEF_S2 strategy continuously releases water for hydropower generation, irrespective of the reservoir storage. This results in emptying of the reservoir earlier and loss of reliability in the MEF objective. It is worth pointing out that the additional hydropower gained by operating the reservoir in this manner is negligible as at lower storage states, the head is also low, resulting in only an incremental gain in hydropower.

**"Insert Figure 5"**

Fig. 5. Hydropower releases (y-axis) versus reservoir storage (x-axis) for strategies PFMEF_S1 and PFnMEF_S2 re-evaluated for observed inflows.

**4.3 Sensitivity of the results to the choice of MEF thresholds**

Sensitivity analysis reveals that the choice of MEF thresholds can have a substantial influence on objective function performance and resultant trade-offs between hydropower generation, water demand deficit, and MEF reliability. We visualise the performance in these objectives across all non-dominated strategies for both PF_MEF and PF_nMEF formulations obtained from re-optimization with MEF thresholds set at 40%, 50%, 60%, 70%, and 80% of MAF (Figure 6). The objective function values are first re-evaluated as a mean across all 100,000 inflow realizations and the resultant values across all non-dominated strategies, after re-sorting across both formulations, are plotted. We find that an increase in MEF requirements for the PF_MEF formulation above the recommended 30% of MAF by Smakhtin and Anputhas (2006) leads to a loss in hydropower generation and demand satisfaction objectives for the PF_MEF formulation. The median hydropower generated decreases from 3662 GWh to 3480 GWh, and the demand deficits increase from 709 $Mm^3$ to 2402 $Mm^3$ when the MEF thresholds are raised from 30% to 80% of MAF for the strategies from the PF_MEF formulation. On the other hand, the median hydropower generated increases from 3457 GWh to 3925 GWh, and the demand deficits increase only slightly from 60 $Mm^3$ to 348 $Mm^3$ when the MEF thresholds are raised from 30% to 80% of MAF for the strategies from the PF_nMEF formulation. An increase in the MEF thresholds above the recommended 30% MAF reduces the MEF reliability for both formulations, however the reduction is worse for the strategies from the PF_nMEF formulation. When the MEF thresholds are increased from 30% to 80% of MAF, the MEF reliability range decreases from 86% - 99% to 66% to 81% for the strategies from the PF_MEF formulation. On the other hand, for the PF_nMEF formulation, as we increase the MEF thresholds from 30% to 80% of MAF, the MEF reliability range decreases from 63% - 80% to 43% to 58%.

We further visualize the trade-offs against the objectives for the 40% and 80% thresholds using a parallel axes plot (Fig. S4). We obtained a total of 32 Pareto-optimal strategies (13 for the PF_MEF, and 19 PF_nMEF formulation) upon re-optimization for the 40% threshold. After re-sorting all strategies for all objective functions, a total of 6 non-dominated strategies for the PF_MEF formulation, and 8 non-dominated strategies for the PF_nMEF formulations were obtained (Fig. S4a). Similarly, when re-optimized for MEF thresholds set at 80% of MAF, we obtained a total of 21 Pareto-optimal strategies (8 for the PF_MEF, and 13 PF_nMEF formulation). A total of 6 non-dominated strategies for the PF_MEF formulation, and 7 non-dominated strategies for the PF_nMEF formulations were obtained on Pareto-sorting (Fig. S4b). Significant losses in hydropower production and demand satisfaction are visible when re-optimized for MEF thresholds adjusted to 40% and 80% of MAF. The losses in these objective performances become the most pronounced when MEF thresholds are 80% of MAF. Therefore, it becomes further difficult to select strategies for PF_MEF formulation with high MEF reliability and acceptable compromises on other objectives on increasing the MEF thresholds beyond 30% of MAF.

Upon re-optimizing both formulations for the MEF thresholds specified by the central body of MoEF&CC, we found the PF_MEF formulation performed much better in the average demand deficit objective with a median value of 319 $Mm^3$. This is still higher than the PF_nMEF strategies that attain a median value of 142 $Mm^3$ but the differences are lesser when compared to more stringent MEF requirements. As the MoEF&CC thresholds require lower MEF releases, enough water is available to meet water demands. However, the hydropower generation is also lower for these new MEF thresholds. This reduction can be attributed to more water being released for demand satisfaction on a priority basis. Thus, there is a significant impact of the choice of MEF thresholds on attainable objective function values and inferred trade-offs.

"Insert Figure 6"

Fig. 6. Boxplots representing objective function ranges of the non-dominated solutions (a) hydropower, (b) average deficit, (c) reliability of non-exceedance of high flows (flood reliability), (d) MEF reliability for different MEF thresholds (MoEF&CC, 40, 50, 60, 70, and 80 percent of MAF) from the PF_MEF and PF_nMEF formulations. The boxplots are color-coded, with green representing the PF_MEF formulation and red representing the PF_nMEF formulation.

## 5. Discussion

Large multipurpose dams store high flows for later use in power generation and water demand satisfaction (Richter and Thomas, 2007). However, these twin goals often come in conflict. A high storage level in the reservoir is favourable for increasing hydropower production that depends both upon the discharge as well as the net head of water. However, maintaining a high head may not be possible as simultaneous release of water for irrigation or domestic use is needed. The conflict between hydropower generation and demand satisfaction objectives has been demonstrated in the Conowingo Dam in the Susquehanna River basin in the USA (Giuliani et al., 2014), the Red River basin in Vietnam (Quinn et al., 2018), and the Tuyamuyun reservoir in the Amu Darya River basin in Uzbekistan (Zhou et al., 2022). In the case of the Tuyamuyun reservoir, for example, an increase in water allocation for irrigation by $1.02 \times 10^9$ m$^3$ results in a 69.9% reduction in water allocation for hydropower generation (Zhou et al., 2022). We also find a significant conflict between these objectives for the case of the NS reservoir. For instance, in the PF_MEF formulation the strategy that performs best in the hydropower objective performs worst in the demand deficit objective and vice-versa. Between these two strategies, increasing the annual hydropower generation by 251 GWh results in a simultaneous increase in annual demand deficit of 974.4 Mm$^3$. These deficits account for 12% of annual demand and the hydropower production value is nearly 7% of maximum possible hydropower generated across PF_MEF strategies.

The number of studies investigating the trade-offs between MEF satisfaction and other socioeconomic objectives have been growing (Feng et al., 2018; Al-Jawad et al., 2019; Quinn et al., 2019; Suwal et al., 2020; Horne et al., 2016; Feng et al., 2018; Marak et al., 2020). These recent studies that consider environmental objectives within a simulation-optimization framework have indicated a possible conflict between demand satisfaction and environmental flow maintenance (Tsai et al., 2015; Dehghanipour et al., 2020). Studies have also shown conflict between environmental flow maintenance and hydropower generation (Li et al., 2016; Feng et al., 2018; Zhang et al., 2019; Jordan et al., 2022). On analysing the reservoir operations downstream of the Lancang River in China, Zhang et al. (2019) showed that a

10.61% reduction in power generation was required to minimally affect ecology downstream. Yan et al. (2021) used the Water Quantity Level (WQL) objective to quantify downstream MEF requirements and reported that an increase hydropower generation of $0.06 \times 10^7$ kWh results in WQL decrease by 0.13. Babel et al. (2012) showed that hydropower can be improved while maintaining MEF requirements in run-of-the-river river plants if the range of MEF thresholds falls in a certain range. In our analysis, we do not observe a significant conflict between MEF reliability and hydropower generation in the case of NS reservoir. For instance, on comparing the best and worst hydropower production strategies for the PF_MEF formulation, an increase in hydropower by 1 GWh leads to a corresponding decrease in MEF reliability by 0.02%. These trade-offs are relatively minor in case of NS dam as the main powerhouse is located at the dam. Thus, hydropower releases contribute to MEF in the main river channel downstream of the dam. Our study further highlights the importance of exploring diverse formulations of water management problems (Quinn et al., 2017; Hadjimichael et al., 2020b). Using the concept of rival framings, Quinn et al. (2017) examined various ways to aggregate objective functions across uncertain realizations and showed that the strategies identified by optimization strongly depended upon the choice of aggregating function. Similarly, Veena et al. (2021) examined multiple problem formulations of a proposed inter-basin water transfer problem involving the NS reservoir to highlight the value of coordinated information exchange between participating basins of the project. Here, we examined how different priority orders of releases in the systems model could affect the optimized strategies attained and inferred multi-objective trade-offs for a multipurpose reservoir. Previous research has emphasized the importance of re-operating dams in a way that mimics the natural flow patterns (Bednarek and Hart, 2005; Richter and Thomas, 2007; Dittmann et al., 2009; Shiau and Wu, 2013). We showed that it is possible to re-operate the NS reservoir by prioritizing MEF while also balancing hydropower generation and demand satisfaction. These operational strategies prioritizing MEF can assist reservoir operators in making informed release decisions.

" While planning and scientific assessments have advanced substantially, successful re-operation of infrastructure to achieve environmental benefits has been more limited (Warner et al., 2014)". Also, understanding the effects of prioritizing environmental flows over other water users has not been explored in previous studies (Wills et al., 2022). Although studies quantifying the ecological benefits of alternative MEF strategies for run-of-the-river plants exist (Basso and Botter, 2012; Kuriqi et al., 2019; Owusu et al., 2023), investigations are limited in multi-purpose reservoirs (Wills et al., 2022). Kuriqi et al. (2019) quantified the trade-offs between multiple objectives of 20 run-of-the-river plants in the Iberian Peninsula and demonstrated that for specific MEF scenarios, a balance between ecological and hydropower objectives can be obtained. On similar lines, we also performed a sensitivity analysis for the NS reservoir by increasing the MEF thresholds (40%, 50%, 60%, 70%, and 80% of MAF). We showed that increasing the MEF thresholds for the PF_MEF formulation beyond the recommended 30% of MAF results in decreased hydropower generation and increased demand deficits. The sensitivity analysis thus sheds insights into the range of MEF thresholds beyond which prioritization of MEF is not likely to yield benefits and lead to conflicts with other objectives. In conclusion, it is important to conduct a case-by-case analysis for testing alternative problem formulations to check practicality of MEF prioritization. Extending this analysis to other reservoirs with varying reservoir metrics such as demand-inflow ratios and impoundment ratios will contribute to developing applicable guidelines MEF for prioritization while reoperating reservoirs. Such analysis also helps in bridging the knowledge gaps between scientists and reservoir operators, thus improving water resources management practices (Richter et al., 2006).

There are a few limitations in this analysis which would need to be addressed to generalize the results obtained here across other multi-purpose reservoirs or their real-world implementation. First, we assumed that water released for hydropower is also used for MEF satisfaction as the main hydropower plant on the NS reservoir is located at the dam site. However, in several instances the power plant will be located away from the dam site, creating a potential conflict between hydropower and MEF releases. Therefore, further problem formulations can be explored to account for this possibility. Second, the implementation of adaptive operating rules that were able to balance all four objectives of the NS reservoir require reservoir operators to follow dynamic rule curves, which are not yet commonly used in India and elsewhere. However, adaptive reservoir operations are also pertinent in the context of shifting water demands and climate change, thus a policy level shift to train reservoir operators to consider such flexible rule curves is needed. Third, previous studies have shown that water demands may decrease during droughts due to farmers changing the cropping patterns or reduction of crop areas irrigated, indicating the importance of incorporating stakeholder inputs and social factors in water management strategies (Venot et al., 2009). In this context, understanding how human-water feedbacks evolve during droughts, and how the water demands change during droughts is important to improve allocation and management of water resources. Finally, it will be necessary to consider the likely impact of a changing climate on inferred success of alternative policies. Climate change will have diverse impacts on water resources, such as alterations in precipitation pattern, increased temperatures, and an increase in the frequency of extreme weather conditions like droughts and floods (Dyer et al., 2014; Jain and Kumar, 2014; Ahn et al., 2018; Zhong et al., 2021). Studies have shown that the anticipated rise in future demands for irrigation and urbanization in major Indian reservoirs will pose a challenging decision-making dilemma for hydropower production and other objectives (Mukheibir et al., 2013; Zeng et al., 2017; IPCC, 2022).

## 6. Conclusions

In this study, two alternative formulations of a multi-purpose reservoir operation decision problem were tested to investigate the effect of the prioritization order of releases on multi-objective performance and inferred trade-offs. By visualising the Pareto approximate solutions for these two alternate formulations we identified the objective trade-offs that can inform policy decisions. We showed that the priority of release decisions plays a significant role in achieving the desired objectives for a multi-objective optimization problem, especially when it comes to attaining required MEF goals. In the formulation that prioritized MEF releases, the use of radial basis functions provided the flexibility needed to balance competing objectives and a policy was discovered that minimally affected water demands while maintaining high performance in MEF satisfaction and hydropower production. Additionally, a sensitivity analysis was conducted to determine the effect of changing the MEF requirements on the resultant multi-objective performance. The results indicated a substantial influence of the choice of MEF thresholds on the value of prioritizing MEF releases; beyond a MEF criteria of 30% MAF prioritizing MEF did not result in well-balanced strategies. Therefore, we aim to develop a systematic framework to quantify the impact of prioritizing MEF releases on other sectoral benefits from multi-purpose reservoirs. This should then enable a consistent data-based approach for MEF prioritization related decisions.

Poff et al. (2009) emphasize the importance of governments setting quantifiable targets for environmental flow management across all rivers and to control other water uses to meet these goals. In this study, we were able to identify the MEF threshold levels for the NS dam re-operation for which we can prioritise for MEF allowing adequate hydropower generation without compromising demand satisfaction. Our sensitivity analysis, incorporating the MEF thresholds based on MoEF&CC guidelines provides evidence supporting the implementation of these thresholds at the reservoir level. Our findings offer valuable information to guide the development of legal frameworks prioritizing MEF. Furthermore, it is crucial to acknowledge that prioritizing MEF in water management goes beyond mere "allocation" and has broader social benefits, which are vital for sustainable water management (Richter, 2010). Lastly, to put the adaptive operational rules from this study into practice, it's important to discourse with a diverse group of stakeholders, including policy makers, engineers, economists, and ecologists.

generation in Central Asia, Journal of Contaminant Hydrology. https://doi.org/10.1016/j.jconhyd.2022.104004.